\documentclass[12pt]{article}
\usepackage{amssymb,amsmath}
\newtheorem{proposicion}{Proposition}[section]

\newtheorem{definicion}[proposicion]{Definition}

\newtheorem{nota}[proposicion]{Remark}

\newtheorem{teorema}[proposicion]{Theorem}
\newenvironment{demostracion}[1]{\paragraph{\sl Proof#1}}{}

\newcommand{\CC}{{\bf C}}

\newcommand{\KK}{{\bf K}}

\newcommand{\NN}{{\bf N}}

\begin{document}

\title{The embedding conjecture for  quasi-ordinary hypersurfaces}
\author{Abdallah Assi\thanks{Universit\'e d'Angers, Math\'ematiques, 49045 Angers cedex 01, France, e-mail:assi@univ-angers.fr}
}
%\date{March 15, 1997}
\date{\mbox{}}
\maketitle
%%\newcommand{\numero}{\refstepcounter{teorema}
%\paragraph{\bf\theteorema}}

%\newtheorem{teorema}{Theorem}[section]
%\newtheorem{teorema}{Theorem}[subsection]

%\renewcommand{\theequation}{\theteorema} (!)
%\newenvironment{formule}{\refstepcounter{teorema}(!)
%\begin{eqnarray}}{\end{eqnarray}}(!)

%$\qed$

%\maketitle

\noindent {\Large {ABSTRACT}}. This paper has two objectives: we first generalize the theory of
Abhyankar-Moh to quasi-ordinary polynomials, then we use the notion of approximate roots and
that of generalized Newton polygons in order to prove the embedding conjecture for this class of
polynomials. This conjecture -made by S.S. Abhyankar and A. Sathaye- says that if a hypersurface of the
affine space is isomorphic  to a coordinate, then it is equivalent
to it.

\section*{Introduction}

\noindent Let ${\bf K}$ be an algebraically closed field of
characteristic zero, and let ${\bf R}={\bf
K}[x_1,\ldots,x_e,y]={\bf K}[\underline{x}][y]$ be the polynomial
ring  $x_1,\ldots,x_e,y$ over ${\KK}$. Let
$f=a_0(x)y^n+a_1(\underline{x})y^{n-1}+\ldots+a_n({\underline{x}})$
be a nonzero  polynomial of ${\bf R}$, and suppose, after a
 possible change of variables that $a_0(\underline{x})=1$. Suppose that
 $e=1$ and that $f$ has  one place at infinity.  Given a nonzero polynomial $g$
of ${\bf R}$,  we define the intersection
 multiplicity of $f$ with $g$, denoted int$(f,g)$, to be the $x_1$-order of the $y$ resultant
 of $f$ and $g$. The set of int$(f,g), g\in {\bf R}$, defines a
 semigroup, denoted $\Gamma(f)$. It follows from the Abhyankar-Moh theory that if
 we denote by  App$(f)=\lbrace g_1,\ldots,g_h\rbrace$ the set of approximate roots of $f$, then
 $n, {\rm int}(f,g_1),\ldots,{\rm
int}(f,g_h)$ generate
 $\Gamma(f)$, and  these polynomials  can be calculated from the equation of $f$ by using the Tschirnhausen
 transform.  As a consequence of this fact, for all $\lambda\in{\bf K}$, App$(f-\lambda)={\rm App}(f)$, and $f-\lambda$ has
 one place at infinity. If furthermore $f$ has no singularities in the affine plane ${\bf K}^2$, then
 it is equivalent to a coordinate, i.e., there is an automorphism $\sigma$ of ${\bf K}^2$ which transforms
 $f$ into a coordinate of ${\bf K}^2$. Let $e\geq 2$ and let $g$ be a nonzero polynomial of ${\bf
R}$, then define the order
 of $g$ with respect to $f$, denoted $\tilde{O}(f,g)$, to be the leading exponent of Res$_y(f,g)$ with respect to
 the diagonal order on ${\bf N}^e$. The set of $\tilde{O}(f,g), g\in{\bf R}[y]$, defines a semigroup, denoted $\Gamma(f)$.
 Let $F(x_1,\ldots,x_e,y)=f(x_1^{-1},\ldots,x_e^{-1},y)\in{\bf K}[x_1^{-1},\ldots,x_e^{-1}][y]\subseteq {\bf
 K}((x_1,\ldots,x_e))[y]$ and call $F$ the meromorphic polynomial associated with $f$.
Suppose that $F$ is irreducible in ${\bf
 K}((x_1,\ldots,x_e))[y]$ and that the discriminant of $F$ is of
the form $x_1^{N_1}.\ldots.x_e^{N_e}.u(x_1,\ldots,x_e)$, where
 $u$ is a unit in ${\bf K}[[\underline{x}]]$ (such a polynomial is called quasi-ordinary polynomial).
 By Abhyankar-Jung Theorem, the roots of $F(x_1,\ldots,x_e,y)=0$ are all in ${\bf K}((x_1^{1\over n},\ldots,x_e^{1\over n}))$, i.e.
 there is a meromorphic series $y(t_1,\ldots,t_e)=\sum_{p\in{\bf Z}^e}c_pt_1^{p_1}.\ldots.t_e^{p_e}\in {\bf K}((t_1,\ldots,t_e))$ such that $F(t_1^n,\ldots,t_e^n,y(t_1,\ldots,t_e))=0$ and any other root of $F(t_1^n,\ldots,t_e^n,y)=0$ is of the form $y(w_1t,\ldots,w_et_e)$, where $w_1,\ldots,w_e$ are $n$th roots of unity in ${\bf K}$. Given a nonzero polynomial $g$ of ${\bf R}$, if $G$ denotes the meromorphic polynomial associated with $g$, we
 define the order of $G$, denoted $O(F,G)$, to be the leading exponent with respect to the lexicographical order
 of the smallest homogeneous component of
$G(t_1^n,\ldots,t_e^n,y(t_1,\ldots,t_e))$.  Clearly $\tilde{O}(F,G)=-O(f,g)$ and the set of $O(F,G)$ defines
a subsemigroup of ${\bf Z}^e$.
The aim of this paper is to use theses notations in order to
generalize the Abhyankar-Moh theory to the set of  polynomials $f$ whose
associated meromorphic polynomials  are irreducible and
quasi-ordinary in ${\bf
 K}((x_1,\ldots,x_e))[y]$.  More precisely we prove the following:

 \medskip

\noindent {\bf Theorem 1.} Let $f$ be a nonzero polynomial of  ${\bf R}$. If  the meromorphic
polynomial  $F$ associated with $f$ is  quasi-ordinary and irreducible in ${\bf K}((x_1,\ldots,x_e))[y]$, then so is for all $F-\lambda, \lambda\in {\bf K}$.

\medskip

\noindent {\bf Theorem 2.} Let $f$ be a nonzero polynomial of ${\bf R}$ and assume that the meromorphic
polynomial $F$ associated with $f$ is  quasi-ordinary  in ${\bf K}((x_1,\ldots,x_e))[y]$.  If $\displaystyle{{\bf R}\over (f)}$ is isomorphic to the
 algebra of a coordinate (i.e. $f$ is isomorphic to a coordinate), then there is a $1\leq k\leq e$ such that $f\in {\bf K}[x_k][y]$. Furthermore, there is an automorphism $\sigma$ of ${\bf K}^{e+1}$ such that $\sigma(f)$ is a
 coordinate of ${\bf K}^{e+1}$ (i.e. $f$ is equivalent to a coordinate).

\medskip

\noindent  Theorem 2. has also the following interpretation: let $f$ be a nonzero polynomial of ${\bf R}$ and suppose
that $f$ is isomorphic to a coordinate, then Abhyankar-Sathaye conjecture  says that $f$
 is equivalent to a coordinate. In particular Theorem 2. gives an affirmative
 answer to the conjecture when the polynomial is quasi-ordinary.

\medskip

 \noindent Our main tools are based on the theory of approximate roots and
 the irreducibility criterion for  polynomials in ${\bf K}((x_1,\ldots,x_e))[y]$.  In a previous
work (see [9]), we proved that given $f$ and $F=f(x_1^{-1},\ldots,x_e^{-1}][y]$, the set of orders of the approximate roots of
$F$ together with the canonical basis of ${\bf Z}^e$ generate the
semigroup of $O(F,G), G\in{\bf
K}[x_1^{-1},\ldots,x_e^{-1}][y]$.  We used this fact and the notion of generalized Newton polygons
in order to give a criterion for the irreducibility of quasi-ordinary polynomials of ${\bf K}((x_1,\ldots,x_e))[y]$.  These
notions are the main tools of the proofs of Theorems 1. and 2.
 Note that the
property of being quasi-ordinary
depends on the choice of coordinates in ${\bf K}^{e+1}$. In this
paper we introduce the notion of almost quasi-ordinary
polynomials. Such polynomials become quasi-ordinary after a
change of variables of ${\bf K}^{e+1}$. In particular
our results remain true for this class of polynomials.

 \noindent The paper is organized as follows: in Section 1 we
 introduce the notion of approximate roots of a polynomial in one
variable  over a commutative ring with unity. In Section 2 we show
how to
 associate a semigroup with an irreducible quasi-ordinary
 polynomial of  ${\bf K}((x_1,\ldots,x_e))[y]$. In Section 3 we adapt the results of Section 2 to the global case. In Section 4 we recall the irreducibility criterion for quasi-ordinary polynomials  in ${\bf K}[x_1^{-1},\ldots,x_e^{-1}][y]$ then we use it in order to prove Theorem 1.  In Section 5.  we use the notion of generalized Newton polygons in order to prove Theorem 2. We also introduce and study the notion of almost quasi-ordinary polynomials.

\section {Approximate roots}

\noindent Let ${\bf S}$ be a commutative ring with unity and let
${\bf S}[y]$ be the ring of polynomials in $y$ with coefficients in
${\bf S}$. Let $f=y^n+a_1y^{n-1}+\ldots+a_n$ be a monic polynomial of
${\bf S}[y]$ of  degree $n > 0$ in $y$. Let $d\in{\NN}$ and suppose
that $d$ divides $n$. Let $g$ be a monic polynomial in ${\bf S}[y]$ of
degree $\displaystyle{n\over d}$ in $y$. There exist unique
polynomials $a_1(y),\ldots,a_d(y)\in{\bf S}[y]$ such that:

$$
 f=g^d+\sum_{i=1}^{d}a_i(y).g^{d-i}
$$

\noindent and for all $1\leq i\leq d$, if we denote by deg$_y$ the
$y$-degree, then deg$_y(a_i)< \displaystyle{n\over d}={\rm
deg}_yg$. The equality above  is called the $g$-adic expansion of
$f$.
\vskip0.2cm

\noindent This construction can be generalized to a sequence of
polynomials. Let to this end $n=d_1 > d_2 >...>d_h$ be a sequence
of integers such that $d_{i+1}$ divides $d_i$ for all $1\leq i\leq
h-1$, and set $e_i=\displaystyle{d_i\over d_{i+1}}$, $1\leq i\leq
h-1$ and $e_h=+\infty$. For all $1\leq i\leq h$, let $g_i$ be a monic polynomial of
${\bf S}[y]$ of degree $\displaystyle{n\over d_i}$ in $y$. Set
$G=(g_1,\ldots,g_h)$ and let $B=\lbrace
(\theta_1,\ldots,\theta_h)\in {\bf
  N}^{h}, 0\leq \theta_i < e_i$ for all $1\leq i\leq h\rbrace$. Then $f$
  can be
 uniquely written in the following form:

$$
 f=\sum_{\underline{\theta}\in
B}a_{\underline{\theta}}.g^{\underline{\theta}}
$$

\noindent where if
$\underline{\theta}=(\theta_1,\ldots,\theta_h)$, then
$g^{\underline{\theta}}=g_1^{\theta_1}.\ldots.g_h^{\theta_h}$ and
$a_{\underline{\theta}}\in{\bf S}$. We call this expansion the
$G$-adic expansion of $f$. We set Supp$_G(f)=\lbrace
\underline{\theta}; a_{\underline{\theta}}\not=0\rbrace$ and we
call it the $G$-support of $f$.

\noindent   Let $f,g$ be as above and let $f=g^d+\sum_{i=1}^{d}a_i(y).g^{d-i}$ be the $g$-adic expansion of $f$. Assume that $d$ is a unit in ${\bf S}$. The Tschirnhausen
transform of $f$ with respect to $g$, denoted $\tau_f(g)$ is
defined by $\tau_f(g)=g+d^{-1}a_1$. Note that $\tau_f(g)=g$ if and
only if $a_1=0$. By [1],  $\tau_f(g)=g$ if and only if
deg$_y(f-g^d) <n-\dfrac{n}{d}$. If one of these equivalent
conditions is verified, then the polynomial $g$ is called a $d$-th
approximate root of $f$. By [1], there exists a unique
 $d$-th approximate root of $f$. We denote it by App$_d(f)$.

\section{The semigroup of a quasi-ordinary polynomial}

\medskip

\noindent Let ${\bf K}$ be an algebraically closed field of
characteristic zero, and let ${\bf A}={\bf K}((x_1,\ldots,x_e))[y]$
(denoted ${\bf K}((\underline{x}))[y]$) be the ring of polynomials in $y$
whose coefficients are meromorphic series  in $x_1,\ldots,x_e$ over ${\KK}$. Let
$F=y^n+a_1(\underline{x})y^{n-1}+\ldots+a_n(\underline{x})$ be
 a nonzero  polynomial of ${\bf A}$ and suppose that the discriminant of $f$ is
 of the form $x_1^{N_1}.\ldots.x_e^{N_e}.u(x_1,\ldots,x_e)$, where
 $N_1,\ldots,N_e\in{\bf Z}$ and $u(\underline{x})$ is a unit
 in ${\bf K}[[\underline{x}]]$. We call $F$  a
 quasi-ordinary polynomial. It follows from Abhyankar-Jung Theorem that there
 exists a  meromorphic  series
 $y(\underline{t})=y(t_1,\ldots,t_e)\in {\bf K}((t_1,\ldots,t_e))$ (denoted ${\bf K}((\underline{t}))$)
 and $m\in {\bf N}$ such that
$F(t_1^m,\ldots,t_e^m,y(\underline{t}))=0$. Furthermore, if $F$ is
an
 irreducible polynomial, then we can take $m=n$, and:

$$
F(t_1^n,\ldots,t_e^n,y)=\prod_{i=1}^n(y-y(w_1^it_1,\ldots,w_e^it_e))
$$

\noindent where  $(w_1^i,\ldots,w_e^i)_{1\leq i \leq n}$ are
distinct elements of $(U_n)^e$, $U_n$ being  the group of $n$th
roots of unity in ${\bf K}$.
\vskip0.2cm
\noindent  Suppose that $F$ is irreducible and let
$y(\underline{t})$ be as above. Write $y(\underline
t)=\sum_{p}c_p{\underline t}^p$ and define the support of $y$ to
be the set $\lbrace p|c_p\not= 0\rbrace$. Obviously the support of
$y(w_1t_1,\ldots,w_et_e)$ does not depend on $w_1,\ldots,w_e\in
U_n$. We denote it by Supp$(F)$ and we call it the support of $F$.
It is well known that there exists a finite sequence of elements
in Supp$(F)$, denoted $m_1,\ldots,m_h$, such that
\vskip0.15cm

i) $m_1 < m_2 <\ldots <m_h$, where $<$ means $<$ coordinate-wise.
\vskip0.15cm
ii) If $c_p\not= 0$, then $p\in (n{\bf Z})^e+\sum_{|m_i|\leq
|p|}m_i{\bf Z}$.
\vskip0.15cm
iii) $m_i\notin  (n{\bf Z})^e+\sum_{j<i}m_j{\bf Z}$ for all
$i=1,\ldots,h$.
\vskip0.15cm
\noindent The set of elements of this sequence is called the set
of characteristic exponents of $F$. We denote by convention
$m_{h+1}=(+\infty,\ldots,+\infty)$. If $e=1$, this set is nothing
but the set of Newton-Puiseux exponents of $F$.

\medskip

\noindent Let $u=\sum_{p}c_p{\underline t}^p$ in ${\bf
K}((\underline t))$ be a nonzero meromorphic  series. We denote by
In$(u)$ the initial form of $u$: if $u=u_d+u_{d+1}+\ldots$ denotes
the decomposition of $u$ into  sum of homogeneous components, then
In$(u)=u_d$. We set $O_t(u)=d$ and we call it the
$\underline{t}$-order of $u$. We denote by exp$(u)$ the greatest
exponent of In$(u)$ with respect to the lexicographical order. We
denote by inco$(u)$ the coefficient $c_{\rm exp}(u)$, and we call
it the initial coefficient of $u$. We set M$(u)={\rm
inco}(u){\underline{t}}^{{\rm exp}(u)}$, and we call it the
initial monomial of $u$.

\vskip0.15cm

\noindent Let $G$ be a nonzero quasi-ordinary element of ${\bf
A}$. The order of $G$ with respect to $F$, denoted $O(F,G)$, is
defined to be  exp($G(t_1^n,\ldots,t_e^n,y(\underline{t})$). Note
that it does not depend on the choice of the root
$y(\underline{t})$ of $F(t_1^n,\ldots,t_e^n,y)=0$. The set
$\lbrace O(F,G)|G\in {\bf A}\rbrace$ defines a subsemigroup of
${\bf Z}^e$. We call it the semigroup associated with $F$ and we
denote it by $\Gamma(F)$.

\medskip

\noindent Let $M(e,e)$ be the unit $(e,e)$ matrix. Let $D_1=n^e$
and for all $1\leq i\leq h$, let $D_{i+1}$ be the gcd of the
$(e,e)$ minors of the matrix $(nM(e,e),{m_1}^T,\ldots,{m_i}^T)$
(where $T$ denotes the transpose of a matrix). Since $m_i\notin
(n{\bf Z})^e+\sum_{j < i}m_j{\bf Z}$ for all $1\leq i\leq h $,
then $D_{i+1}< D_{i}$. We define the sequence $(e_i)_{1\leq i \leq
h}$ to be $\displaystyle{e_i={D_i\over D_{i+1}}}$ for all $1\leq
i\leq h$.

\noindent Let $M_0=({n\bf Z})^e$ and let $M_i= (n{\bf
Z})^e+\sum_{j=1}^im_j{\bf Z}$ for all $1\leq i\leq h$. Then $e_i$
is the index of the lattice $M_{i-1}$ in $M_i$, and
$n=e_1.\ldots.e_h$, in particular $D_{h+1}=n^{e-1}$. We set
$d_i=\displaystyle{{D_i}\over {D_{h+1}}}$ for all $1\leq i\leq
h+1$. In particular $d_1=n$ and $d_{h+1}=1$. The sequence
$(d_1,d_2,\ldots,d_{h+1})$ is called the gcd-sequence associated
with $f$. We also define the sequence $(r_k)_{1\leq k \leq  h}$ by
$r_1=m_1$ and $r_{k+1}=e_kr_k+m_{k+1}-m_k$ for all $1\leq k \leq
h-1$.

\vskip0.2cm

 \noindent Let
 $\phi(\underline t)=(t_1^p,\ldots,t_e^p,Y(\underline{t}))$ and
 $\psi(\underline t)=(t_1^q,\ldots,t_e^q,Z(\underline{t}))$
be two nonzero elements of ${\bf K}((\underline t))^{e+1}$. We
define the contact between $\phi$ and $\psi$ to be the element
${\dfrac{1}{pq}}{\rm
exp}(Y(t_1^q,\ldots,t_e^q)-Z(t_1^p,\ldots,t_e^p))$. We denote it
by $c(\phi,\psi)$.

\noindent We define the contact between $F$ and $\phi$, denoted
$c(f,\phi)$, to be the maximal element in the set of contacts
of $\phi$ with the roots $(t_1^n,\ldots,t_e^n,y(\underline{t})$ of $F(x_1,\ldots,x_e,y)=0$.

\noindent Let
$G=y^m+b_1(\underline{x})y^{m-1}+\ldots+b_m(\underline{x})$ be a
nonzero polynomial of ${\bf A}$. Suppose that $G$ is an
irreducible quasi-ordinary polynomial and let $\psi(\underline
t)=(t_1^m,\ldots,t_e^m,Z(\underline{t}))$ be a root of
$G(x_1,\ldots,x_e,y)=0$. We define the contact between $F$ and
$G$, denoted $c(F,G)$, to be the contact between $F$ and $\psi$,
 and we recall that this definition
does not depend on the choice of the root $\psi$ of $G$. Note that
if $F.G$ is a quasi-ordinary polynomial, then ${\rm
In}(F(\psi(\underline{t}))=M(F(\psi(\underline{t}))$.

\vskip0.15cm

\noindent With these notations we have the following proposition:
\vskip0.15cm

\begin{proposicion}{\rm  Let $G=y^m+b_1(\underline{x})y^{m-1}+\ldots+b_m(\underline{x})$ be an
 irreducible quasi-ordinary polynomial of ${\bf
A}$ and suppose that $F.G$ is a quasi-ordinary polynomial. Let
$(D'_j)_{1\leq j\leq h'+1}$ (resp. $(d'_j)_{1\leq j\leq h'+1}$,
$(m'_j)_{1\leq j\leq h'}$)  be the set of characteristic sequences
associated with $G$. If $c$ denotes the contact $c(F,G)$, then we
have the following:
\vskip0.15cm
i) If for all $1\leq 1\leq h$, $nc\notin M_q$, then $\tilde{O}(F,G)=n.m.c$.
\vskip0.15cm
ii) Otherwise, let $1\leq q\leq h$ be the smallest integer
such that $nc\in M_{q}$, then
$\tilde{O}(f,g)=(r_qd_q+(nc-m_q)d_{q+1}).{\dfrac{m}{n}}$.
\vskip0.15cm

iii) If $nc\in M_{q}-M_{q-1}$  and $nc\not=m_q$, then
$\dfrac{n}{d_{q+1}} |m$.}\end{proposicion}

\begin{demostracion}{.}{\rm i) and ii) are obvious. To prove iii) let $\phi=(t_1^n,\ldots,t_e^m,Y(\underline{t}))$
(resp. $\psi=(t_1^m,\ldots,t_e^m,Z(\underline{t}))$) be a root of
$F(\underline{x},y)=0$ (resp. $G(\underline{x},y)=0$) and  remark
that if $nc\in M_{q}-M_{q-1}$ and $nc\not= m_q$ then the exponents of
$Z({t_1}^n,...,{t_e}^n)$ coincide with those of
$Y({t_1}^m,...,{t_e}^m)$ till at least $m_q.m$. Write
$Y(\underline{t})=\sum_ic_i{\underline{t}}^i$ and
$Z(\underline{t})=\sum_jc'_j\underline{t}^j$, then for all $i\in
M_{q+1}$ in Supp$(Y)$, there exists $j\in {\rm Supp}(Z)$ such that
$i.m=j.n$. But the gcd of minors of the matrix
$(m.nM(e,e),t_{m.m_1},\ldots,t_{m.m_q})$ is $m^e.D_{q+1}$, and the
gcd of minors of the matrix
$(m.nM(e,e),t_{n.m'_1},\ldots,t_{n.m'_q})$ is
 $n^e.D'_{q+1}$. Thus
$m^e.D_{q+1}=n^e.D'_{q+1}$, in particular
$m^e.n^{e-1}d_{q+1}=n^e.m^{e-1}.d'_{q+1}$. This implies that
$m=\displaystyle{n\over d_{q+1}}.d'_{q+1}$, which proves our
assertion.}
\end{demostracion}

\vskip0.15cm

\noindent  In [9], it was proved that if $F\in{\bf K}[[x_1,\ldots,x_e]]$, then for
all $k=1,\ldots,h$, $O(F,{\rm App}_{d_k}(F))=r_k$. The same proof works in the
general case, more precisely we have the following:

\begin{teorema}{\rm Let the notations be as above, and let $d_1,\ldots,d_{h}, d_{h+1}=1$ be the
gcd-sequence of $F$. Then for all $1\leq k\leq h$ we have:

\vskip0.15cm

i)  $c(F,{\rm App}_{d_k}(F)=m_k$
\vskip0.15cm

ii) $O(F,{\rm App}_{d_k}(F))=r_k$.}
\end{teorema}

\section{The polynomial case}

\medskip

\noindent Let the notations be as in the introduction, namely ${\bf R}$ denotes the ring of polynomials
in $x_1,\ldots,x_e,y$ with coefficients in ${\bf K}$, and  $f=y^n+a_1(\underline{x})y^{n-1}+\ldots+a_n(\underline{x})$ is
a nonzero element of ${\bf R}$. Let $F=f(x_1^{-1},\ldots,x_e^{-1},y)\in {\bf K}[x_1^{-1},\ldots,x_e^{-1}][y]$ be the meromorphic polynomial associated with $f$, and suppose that $F$ is an irreducible, quasi-ordinary polynomial in ${\bf K}((x_1,\ldots,x_e))[y]$. Let $g$ be a nonzero element of ${\bf R}$ and let $G$ be the meromorphic polynomial associated with $g$.  Clearly Res$(F,G,y)$ is the meromorphic polynomial associated with Res$(f,g,y)$. Let $\tilde{O}(f,g)$ denotes the leadins exponent of Res$(f,g,y)$ with respect to the diagonal order. The set of $\tilde{O}(f,g), g\in {\bf R}-\lbrace 0\rbrace$ is a semigroup. We call it the semigroup of $f$, and we denote it by $\tilde{\Gamma}(f)$.  Similarly we define the semigroup of $F$, denoted $\Gamma(F)$,  to be the set of $O(F,G), G\in {\bf K}[x_1^{-1},\ldots,x_e^{-1}][y]$. Let $(d_k)_{1\leq k\leq h}, (r_k)_{1\leq k\leq h}$, and $(m_k)_{1\leq k\leq h}$ be the set of characteristic sequences associated with $F$ as in Section 2. We have the following theorem:

\begin{teorema}{\rm With the notations be as above, we have the following:
\vskip0.15cm
i)  $\tilde{\Gamma}(f)$ is generated by  $-r_0^1=(n,0,\ldots,0),\ldots,-r_0^e=(0,\ldots,0,n), -r_1,\ldots,-r_h$.
\vskip0.15cm
ii)  For all $1\leq k\leq h$, App$_{d_k}(F)$ is the meromorphic polynomial associated with App$_{d_k}(f)$. In particular
$\tilde{O}(f,{\rm App}_{d_k}(f))=-r_k$.
\vskip0.15cm
iii) For all $1\leq k\leq h$, App$_{d_k}(F)$ is an irreducible quasi-ordinary polynomial in ${\bf K}((x_1,\ldots,x_e))[y]$, of degree $\displaystyle{n\over d_k}$ in $y$, and $\Gamma({\rm App}_{d_k}(F))$ is generated by the caconical basis of $\displaystyle({n\over d_k}{\bf Z})^e$
and $\displaystyle{{r_1\over d_{k-1}},\ldots,{r_k\over d_{k-1}}}$
 }\end{teorema}

\begin{demostracion}{.} {\rm For all $g\in{\bf R}$, we have $\tilde{O}(f,g)=-O(F,G)$, where $G$ is the meromorphic polynomial associated with $g$. This with Theorem 2.2. imply our results.}
\end{demostracion}

\section{Generalized Newton polygons}

\medskip
\noindent In the following we shall recall the notion of generalized Newton polygons of an irreducible quasi-ordinary polynomial.
This notion has been used in order to give a criterion for a quasi-ordinary polynomial of ${\bf K}[[x_1,\ldots,x_e]][y]$ to be irreducible. This criterion works also in the general case, namely for quasi-ordinary polynomials in  ${\bf K}((x_1,\ldots,x_e))[y]$.
More precisely  let $n\in {\bf N}$ and let
$-{\underline{r}}_0=(-r_0^1,\ldots,-r_0^e)$ be the canonical basis of
$(n{\bf Z})^e$. Let $r_1>\ldots >r_h$ be a sequence of elements of
${\bf -N}^e$, where $>$ means $>$ coordinate-wise. Set $D_1=n^e$
and for all $1\leq k\leq h$, let $D_{k+1}$ be the GCD of the
$(e,e)$ minors of the $(e,e+k)$ matrix
$(n.I(e,e),(r_1)^T,\ldots,(r_k)^T)$. Suppose that $n^{e-1}$
divides $D_{k}$ for all $1\leq k\leq h+1$ and that
$D_{h+1}=n^{e-1}$, and also that $D_1
> D_1 >\ldots > D_{h+1}$, in such a way that if we set $d_1=n$ and
$\displaystyle{d_k={D_k\over n^{e-1}}}$ for all $2\leq k \leq h$,
then $d_1=n > d_2 >\ldots > d_{h+1}=1$.

\medskip

\noindent  For all $1\leq k\leq h$, let $G_k\in{\bf K}[x_1^{-1},\ldots,x_e^{-1}][y]$ be a monic
polynomial of degree $\displaystyle{n\over d_k}$ in $y$ and set
$\underline{G}=(G_1,\ldots,G_h)$. Let $F$ be a nonzero polynomial of ${\bf
K}[x_1^{-1},\ldots,x_e^{-1}]][y]$ and let:

$$
F=\sum_{\underline{\theta}\in
B(\underline{G})}c_{\underline{\theta}}(\underline{x})G_1^{\theta_1}.\ldots.G_h^{\theta_h}
$$

\noindent where $B(\underline{G})=\lbrace \underline{\theta}=
(\theta_1,\ldots,\theta_h); \forall 1\leq i\leq h-1, 0\leq
\theta_i < e_i=\dfrac{d_i}{d_{i+1}}$ and $\theta_{h} <
+\infty\rbrace$, be the $\underline{G}$-adic expansion of $F$. Let
Supp$_{\underline{G}}(F)=\lbrace \underline{\theta}\in B(\underline{G});
c_{\underline{\theta}}\not= 0\rbrace$. If $\theta\in {\rm
Supp}_{\underline{G}}(F)$ and $\underline{\gamma}={\rm
exp}(c_{\underline{\theta}}(\underline{x}))$, we shall associate
with the monomial
$c_{\underline{\theta}}(\underline{x})g_1^{\theta_1}.\ldots.g_{h}^{\theta_{h}}$
the $e$-uplet

$$
<((\underline{\gamma},\underline{\theta}),({\underline{r}}_0,\underline{r}))>=\sum_{i=1}^e\gamma_i.r_0^i+\sum_{j=1}^{h}\theta_j.r_j
$$

\noindent There is a unique $\underline{\theta}^0\in {\rm
Supp}_{\underline{G}}(F)$ such that if ${\underline{\gamma}}^0={\rm
exp}(c_{{\underline{\theta}}^0}(\underline{x}))$, then:

$$
<(({\underline{\gamma}}^0,{\underline{\theta}}^0),({\underline{r}}_0,\underline{r}))>=
{\rm inf}\lbrace
<((\gamma,\underline{\theta}),({\underline{r}}_0,\underline{r}))>,
\underline{\theta}\in {\rm Supp}_{\underline{G}}(F)\rbrace
$$

\noindent  We set

$$
{\rm
fO}(\underline{r},{\underline{G}},F)=<(({\underline{\gamma}}^0,{\underline{\theta}}^0),({\underline{r}}_0,\underline{r}))>
$$

\noindent and we call it the formal order of $F$ with respect to
$(\underline{r},\underline{G})$.

% We also set:

%$$
%M_G(F)=M(c_{{\underline{\theta}}_0}).g_1^{\theta_1^0}.\ldots.g_h^{\theta_h^0}
%$$

%\noindent and we call it the initial monomial of $F$ with respect
%to $(\underline{r},\underline{G})$.

\medskip

 \noindent Suppose that $F$ is monic in $y$, then write
$F=y^n+a_1(\underline{x})y^{n-1}+\ldots+a_n(\underline{x})$. Suppose that $F$ is a
quasi-ordinary polynomial of ${\bf K}((x_1,\ldots,x_e))[y]$ and let $d\in {\bf N}$ be a divisor
of $n$. Let $G$ be a monic
polynomial of ${\bf K}[x_1^{-1},\ldots,x_e^{-1}][y]$ of degree
$\displaystyle{{n\over d}}$ in $y$ and let:

$$
F=G^{d}+a_1(\underline{x},y)G^{d-1}+\ldots+a_{d}(\underline{x},y)
$$

\noindent be the $G$-adic expansion of $F$. We associate with $F$
the set of points:

$$
\lbrace ({\rm fO}(\underline{r},\underline{G},a_k),(d-k){\rm
fO}(\underline{r},\underline{G},G)), k=0,\ldots,d\rbrace \subseteq {\bf
N}^{e}\times {\bf N}^e
$$

\noindent We denote this set by GNP$(F,\underline{r},\underline{G},G)$ and we
call it the generalized Newton polygon of $F$ with respect to
$(\underline{r},\underline{G},G)$.
%Note that if $e=1$ and $F$ is an
%irreducible polynomial of ${\bf K}[[x]][y]$, then the above set is
%equivalent to the usual Newton polygon of $f$.

\begin{definicion}{\rm We say that $F$ is straight with respect to $(\underline{r},\underline{G},G)$ if the following holds:

\vskip0.15cm

i) fO$((\underline{r},\underline{G},a_d)=d.{\rm fO}((\underline{r},\underline{G},G))$.

\vskip0.15cm

ii) For all $1\leq k\leq h-1$, fO$(\underline{r},\underline{G},a_k)\geq k.{\rm
fO}((\underline{r},\underline{G},G))$, where $\geq$ mean $\geq$
coordinate-wise.
\vskip0.15cm

\noindent We say that $F$ is strictly straight with respect to
$(\underline{r},\underline{G},G)$ if the inequality in ii) is a strict
inequality.}
\end{definicion}

 \noindent {\bf The criterion} Let $f=y^n+a_1(x)y^{n-1}+\ldots+a_n(x)$  be a nonzero element of
 ${\bf K}[x_1,\ldots,x_e][y]$ and assume, possibly after a change of
variables, that $a_1(\underline{x})=0$.  Let $F\in {\bf K}[x_1^{-1},\ldots,x_e^{-1}][y]$ be the meromorphic polynomial associated with $f$. Let $-{\underline{r}}_0=(-r^1_0,\ldots,-r^e_0)$ be the canonical basis of
$(n{\bf Z})^e$ and let $d_1=n$. Let $G_1=y$ be the $d_1$-th
approximate root of $F$ and set $m_1=r_1={\rm
exp}(a_n(\underline{x}))$. Let $D_2$ be the gcd of the $(e,e)$
minors of the $(e,e+1)$ matrix $(n.I(e,e),{m_1}^T)$. Let
$d_2=\displaystyle{D_2\over n^{e-1}}$ and let $G_2$ be the
$d_2$-th approximate root of $F$ and set
$e_2=\displaystyle{{d_1\over d_2}={n\over d_2}}$.... Suppose that
we constructed $(r_1,\ldots,r_{k-1})$, $(m_1,\ldots,m_{k-1})$, and
$(d_1,\ldots,d_{k})$, then let $G_k$ be the $d_k$-th approximate
root of $F$ and let

$$
F=G_k^{d_k}+\beta_2^kG_k^{d_k-2}+\ldots+\beta_{d_k}^k
$$

\noindent be the $G_k$-adic expansion of $f$. Then $r_k={\rm
fO}(\underline{r}^k,{\underline G}^k,\beta_{d_k}^k)$, where
$\underline{r}^k=(\dfrac{r^1_0}{d_k},\ldots,\dfrac{r^e_0}{d_k},
\dfrac{r_1}{d_k},\ldots,\dfrac{r_{k-1}}{d_k})$ and
${\underline{G}}^k=(G_1,\ldots,G_{k-1})$. With these notations we have the
following:

\begin{teorema}{\rm The polynomial $F$ is an irreducible quasi-ordinary polynomial
in ${\bf K}((x_1,\ldots,x_e))[y]$ if and only if the following holds:

\vskip0.15cm

i) There is an integer $h$ such that $d_{h+1}=1$.

\vskip0.15cm

ii) For all $1\leq k\leq h-1, r_kd_k < r_{k+1}d_{k+1}$, where $<$
means $<$ coordinate-wise.
\vskip0.15cm

iii) For all $2\leq k\leq h+1$, $G_k$ is strictly straight with
respect to $(\underline{r}^k,{\underline{G}}^k,G_{k-1})$.}
\end{teorema}

\vskip0.15cm

\noindent As a corollary we get our first main Theorem:
\vskip0.15cm
\begin{teorema}{\rm Let $f$ and $F$ ba as above. If $F$ is an irreducible quasi-ordinay polynomial of ${\bf K}((x_1,\ldots,x_e))[y]$ then  for all $\lambda\in {\bf K}$ we have the following:
\vskip0.15cm
i) $\Gamma(F-\lambda)=\Gamma(F)$.
\vskip0.15cm
ii) App$_{d_k}(F-\lambda)={\rm App}_{d_k}(F)$ for all $k=1,\ldots,h$.
\vskip0.15cm
\noindent Furthermore, $F-\lambda$ is an irreducible quasi-ordinary polynomial of ${\bf K}((x_1,\ldots,x_e))[y]$.}
\end{teorema}
\begin{demostracion}{.} Note that if $F$ is a quasi-ordinary polynomial, then so is for $F-\lambda$ for all $\lambda\in {\bf K}$. On the other hand, conditions i) and ii) are obvious, then a direct application of the criterion above shows that $F-\lambda$ is irreducible in ${\bf K}((x_1,\ldots,x_e))[y]$.
\end{demostracion}

\section{Quasi-ordinary polynomials isomorphic to a coordinate}
\medskip
\noindent  Let $f=y^n+a_1(x)y^{n-1}+\ldots+a_n(x)$  be a nonzero element of
 ${\bf K}[x_1,\ldots,x_e][y]$ and assume, possibly after a change of
variables, that $a_1(\underline{x})=0$.  Let $F\in {\bf K}[x_1^{-1},\ldots,x_e^{-1}][y]$ be the meromorphic
polynomial associated with $f$ and assume that $F$
is an irreducible  quasi-ordinary polynomial of ${\bf K}((x_1,\ldots,x_e))[y]$. Let $r_0^1=(-n,0,\ldots,0),\ldots,r_0^e=(0,\ldots,0,-n)$
and let $r_1,\ldots,r_h$ be such that $\Gamma(F)=<r_0^1,\ldots,r_0^e,r_1,\ldots,r_h>$.
Let $d\in {\bf N}$ be a divisor of $n$ and let $g$ be a monic
polynomial of ${\bf K}[x_1,\ldots,x_e][y]$ of degree
$\displaystyle{{n\over d}}$ in $y$. Let:

$$
f=g^{d}+a_1(\underline{x},y)g^{d-1}+\ldots+a_{d}(\underline{x},y)
$$

\noindent be the $g$-adic expansion of $f$ and consider the set of points:

$$
\lbrace (\tilde{O}(f,a_k),(d-k)
\tilde{O}(f,g)), k=0,\ldots,d\rbrace \subseteq {\bf
N}^{e}\times {\bf N}^e
$$

\noindent Similarly with Definition 4.1., we say that  $f$ is straight
(resp. stricly straight) with respect to $g$ if the following holds:

\medskip

i) $\tilde{O}(f,a_d)=d.{\tilde{O}}(f,g)$.

\medskip

ii) For all $1\leq i\leq d-1$,  $\tilde{O}(f,a_i) \leq i.(-r_i)$ (resp.  $\tilde{O}(f,a_i)  <  i.(-r_i)$
\medskip

\noindent  Let $g_k={\rm App}_{d_k}(f)$ for all $k=1,\ldots,h+1$ (in particular $g_1=y$ and $g_{h+1}=f$). It results from Section 4. that  $g_{k+1}$ is strictly straight with respect to $g_k$ for all $1\leq k\leq h$. More precisely,
let $1\leq k\leq h$ and set $\displaystyle{e_k={d_k\over d_{k+1}}}$. If:

$$
g_{k+1}=g_k^{e_k}+a^k_2(\underline{x},y)g_k^{e_k-2}+\ldots+a^k_{e_k}(\underline{x},y)
$$

\noindent denotes the $g_k$-adic expansion of $g_{k+1}$, then we have:

\medskip

i) $\tilde{O}(g_{k+1},a^k_{e_k})=\displaystyle{e_k.\tilde{O}(g_{k+1},g_k)=e_k.(-{r_{k+1}\over d_k}})$.

\medskip
ii) For all $2\leq i\leq e_k-1$,  $\tilde{O}(f,a^k_i)  <  i.\displaystyle{(-{r_{k+1}\over d_k})}$.

\medskip

\noindent  Assume that $f$ is equivalent to a coordinate, then
${\bf R}/(f)\simeq {\bf K}[t_1,\ldots,t_e]$, and the meromorphic
polynomial $F$ associated with $f$ is irreducible in ${\bf
K}((x_1,\ldots,x_e))[y]$. It follows that $\tilde{\Gamma}(f)$
contains the canonical basis of ${\bf Z}^e$. If $n=1$, then
$f=z+g(x_1,\ldots,x_e)$ which is equivalent to a coordinate.
Suppose that $n > 1$: With the notations of Section 3., we have
$-r_1d_1
> \ldots
> -r_hd_h$, in particular there is $1\leq k\leq e$ such that $-r_h=(0,\ldots,0,1,0,\ldots,0)$ is the $k$th element
 of the canonical basis of ${\bf Z}^e$. Let:

$$
f=g_h^{d_h}+a^h_2(\underline{x},y)g_h^{d_h-2}+\ldots+a^h_{d_h}(\underline{x},y)
$$

\noindent be the $g_h$-adic expansion of $f$. The above conditions give the following:

\vskip0.15cm

i) $\tilde{O}(f,a^h_{d_h})=d_h.\tilde{O}(f,g_h)=-r_hd_h=(0,\ldots,0,d_h,0,\ldots,0)$.

\vskip0.15cm

ii) For all $2\leq i\leq d_h-1$,  $\tilde{O}(f,a^h_i)  <  -i.r_h=(0,\ldots,i,0,\ldots,0)$.

\vskip0.15cm

\noindent Since $\tilde{O}(f,a^h_i)\in <-r_0^1,\ldots,-r_0^e,-r_1,\ldots,-r_{h-1}>$, and $i < d_h$, condition ii) implies that
$\tilde{O}(f,a^h_i)=0$, i.e. $a^h_i\in {\bf K}$ for all $2\leq i\leq d_h-1$. On the other hand,  $\tilde{O}(f,a^h_i)\in <-r_0^1,\ldots,-r_0^e,-r_1,\ldots,-r_{h-1}>$ implies that there exist $b_0^1,\ldots,b_0^e,b_1,\ldots,b_h\in {\bf N}$ such that:

$$
(-r_h)d_h=(0,\ldots,0,d_h,0,\ldots,0)=\tilde{O}(f,a^h_i)=b_0^1(-r_0^1)+\ldots+b_0^e(-r_0^e)+b_1(-r_1)+\ldots+b_{h-1}(-r_{h-1})
$$

\noindent Since $a^h_{d_h}\not= 0$, this is possible only if $b_0^1=\ldots=b_0^e=b_1=\ldots=b_{h-2}=0, b_{h-1}=1$, and $r_{h-1}=(0,\ldots,0,d_h,0,\ldots,0)$. This implies in particular that $a^h_{d_h}=c_h.g_{h-1}$, where $c_h\in {\bf K}^*$. Finally:

\medskip

1.  $f=g_h^{d_h}+a^h_2(\underline{x},y)g_h^{d_h-2}+\ldots+c_h.g_{h-1}$ and $(a^h_2,\ldots,a^h_{d_h-1},c_h)\in {\bf K}^{d_h-2}\times {\bf K}^*$.

\medskip

2. $-r_h=(0,\ldots,0,1,0\ldots,0)$ and $-r_{h-1}=(0,\ldots,0,d_h,0,\ldots,0)$ (where $1$ and $d_h$ are at the $k$th place).

\medskip

\noindent Condition 2. implies that $\displaystyle{-{r_{h-1}\over d_h}}=-r_h$ is the $k$th element of the canonical basis of
${\bf Z}^e$. But

$$
\displaystyle{\tilde{\Gamma}(g_h)=<-{r_0^1\over d_h},\ldots,-{r_0^e\over d_h},-{r_1\over d_h},\ldots,-{r_{h-1}\over d_h}>}
$$

\noindent in particular the same argument as above applied to $g_h$ instead of $f$ implies that:

\vskip0.15cm

1.  $g_h=g_{h-1}^{e_h}+a^{h-1}_2(\underline{x},y)g_{h-1}^{e_h-2}+\ldots+c_{h-1}.g_{h-2}$ and $(a^{h-1}_2,\ldots,a^{h-1}_{e_h-1},c_{h-1})\in {\bf K}^{e_h-2}\times {\bf K}^*$.

\vskip0.15cm

2. $-r_{h-1}=(0,\ldots,0,d_h,0\ldots,0)$ and $-r_{h-2}=(0,\ldots,0,d_{h-1},0,\ldots,0)$, where $d_h$ and $d_{h-1}$ are at the $k$th place.

\medskip

\noindent Now an easy argument of  induction  implies that for all $1\leq i\leq h$ we have:
\medskip

1.  $-r_i=(0,\ldots,0,d_{i+1},0,\ldots,0)$, where $d_{i+1}$ is at the $k$th place.

\medskip

2.  $g_{i+1}=g_i^{e_i}+a^i_2(\underline{x},y)g_i^{e_i-2}+\ldots+c_ig_{i-1}$ and $(a^{i}_2,\ldots,a^{i}_{e_i-1},c_{i})\in {\bf K}^{e_i-2}\times {\bf K}^*$

\medskip

\noindent in particular  $g_1=y$ and $g_2=y^{e_1}+a^1_2y^{e_1-2}+\ldots+a^1_1y+c_0x_k$, where $(a^1_2,\ldots,a^1_{e_1-1},c_0)\in {\bf K}^{e_1-2}\times {\bf K}^*$. It also follows that $g_i\in {\bf K}[y,x_k]$ for all
$i=1,\ldots,h+1$. Furthermore, ${\bf K}[f,g_h]=\ldots={\bf K}[y,x_k]$. This result can be announced as follows:

\begin{teorema}{\rm Let $f=y^n+a_1y^{n-1}+\ldots+a_n$ be a nonzero element of ${\bf R}$. Suppose, after a possible change of variables, that $a_1=0$ and let $F$ be the meromorphic polynomial associated with $f$. If $f$ is equivalent to a coordinate and if $F$ is a quasi-ordinary polynomial  in ${\bf K}((x_1,\ldots,x_e))[y]$, then we have the following:

\medskip

i) There is $1\leq k\leq e$ such that $f\in {\bf K}[x_k,y]$.

\medskip

ii) There is an automorphism that transforms $f$ into a coordinate of ${\bf K}[x_k,y]$, hence of ${\bf K}[x_1,\ldots,x_e][y]$.}
\end{teorema}

\vskip0.15cm

\begin{nota}{\rm   Since ${\bf K}[x,y,z]={\bf K}[w=x+y,y,z]={\bf K}[z^3-w,y,z]$, then the polynomial $f=z^3-x-y$ is equivalent to a coordinate. However, $f$ is not quasi-ordinary with respect to any of the variables $x,y,z$. Now the change of variables $w=x+y,y,z$ transforms $f$ into $z^3-w$, which is quasi-ordinary. This example suggests to introduce a new class of polynomials, more precisely, given a nonzero polynomial $f=y^n+a_1(\underline{x})y^{n-1}+\ldots+a_n(\underline{x})$ of ${\bf R}$, we say that $f$ is almost quasi-ordinary (a.q.o. for short) if there is a change of variables of ${\bf K}[x_1,\ldots,x_e]$ such that $f$ becomes quasi-ordinary in the new coordinates, say $w_1,\ldots,w_e,y$. A quasi-ordinary polynomial is clearly almost quasi-ordinary, furthermore, given an almost quasi-ordinary polynomial, if $f$ is isomorphic to a coordinate, then $f$ is equivalent to it. This covers a larger set of polynomials. Note that the property of being quas-ordinary is effective, i.e. we can easily decide from the equation if a nonzero polynomial $f$ of ${\bf R}$ is quasi-ordinary. This is not the case  for the almost quasi-ordinary  property, since the structure of automorphisms of ${\bf K}^e, e\geq 3$ is not well understood. Since automorphisms of ${\bf K}^2$ are well known it is natural to consider the case of surfaces, more precisely  we will consider and  answer the following question: Given a nonzero polynomial $D(X,Y)\in {\bf K}[X,Y]$, can we decide effectively if there is a change of variables $X=d_1(w_1),Y=d_2(w_2)$ such that in the new coordinates $w_1,w_2,
D(w_1^{-1},w_2^{-1})=w_1^{N_1}w_2^{N_2}(1+u(w_1,w_2)$, where $N_1, N_2 \leq 0$ and $u(0,0)=0$ (in this case we say that $D$ has the q.o. property)? the main tool here is the Newton polygon of $D(X,Y)$: write $D(X,Y)=\sum_{ij}a_{ij}X^iY^j$ and  define the Newton polygon of $D$, denoted $N(D)$, to be the convex hull of the set $\lbrace (0,0)\rbrace \bigcup \lbrace (i,j), a_{ij}\not=0\rbrace$. Let $E$ be the set of
edges of $P(D)$ with negative slopes. We have the following:
\vskip0.05cm

\noindent 1. If $E=\emptyset$ then $D(X,Y)=a_{i_0j_0}.X^{i_0}Y^{j_0}+\sum_{i <i_0, j <j_0}a_{ij}X^iY^j$ and $D$ has the q.o. property.
\vskip0.05cm
\noindent 2. Suppose that  $E\not= \emptyset$ and let $E=\lbrace E_1,\ldots,E_s\rbrace$.
\vskip0.05cm

2.1. If  s$\geq 2$ then $D$ does not have  the q.o. property.
\vskip0.05cm

2.2.  Suppose that $s=1$ and let $P(E_1)=\sum_{(i,j)\in E_1,a_{ij}\not=0}a_{ij}X^iY^j$. Let $P(E_1)=P_1^{p_1}.\ldots.P_r^{p_r}$ be the decomposition of $P(E_1)$ into irreducible quasi-homogeneous polynomials.
\vskip0.05cm

\quad 2.2.1. If $r\geq 3$ then $D$does not have the q.o. property.
\vskip0.05cm
\quad 2.2.2. If $r=2$ and ${\bf K}[P_1,P_2]\not={\bf K}[X,Y]$ then  $D$ does not have the q.o. property.
\vskip0.05cm
\quad 2.2.3 If $r=1$ and $P_1$ is not equivalent to a coordinate, then $D$ does not have the q.o. property.
\vskip0.05cm
\quad 2.2.4. If $r=2$ and  ${\bf K}[P_1,P_2]={\bf K}[X,Y]$, then wo apply the change of variables $(X,Y)\longmapsto (P_1,P_2)$.
\vskip0.05cm
\quad 2.2.5. If $r=1$ and $P_1$ is equivalent to a coordinate, let $Q$ be such that ${\bf K}[P_1,Q]={\bf K}[X,Y]$,  then we apply the change of variables $(X,Y)\longmapsto (P_1,Q)$.
\vskip0.05cm
\noindent We easily verify that the algorithm above answers our question.}\end{nota}

 \begin{nota}{\rm Let in general $f$ be a polynomial in $s$ variables $X_1,\ldots,X_s$ over ${\bf K}$. Abhyankar-Sathaye conjecture says the following: if $f$ is isomorphic to a coordinate, then $f$ is equivalent to it. This
 conjecture is known to be true for $s=2$ -by Abhyankar-Moh results-  but it is still open for $s\geq 3$.  Theorem 5.1. and
 Remark 5.2. give an affirmative answer to the conjecture
 for almost quasi-ordinay polynomials.}
\end{nota}
%\bigskip

%\begin{thebibliography}{11}\end{thebibliography}

%\centerline {\bf References}

\end{document}